\documentclass[12pt]{article}
\usepackage{theorem,amsfonts}
\newtheorem{theorem}{Theorem}[section]
\newtheorem{proposition}{Proposition}[section]
\newtheorem{lemma}{Lemma}[section]
\newtheorem{corollary}{Corollary}[section]
\theorembodyfont{\upshape}
\newtheorem{definition}{Definition}
\newtheorem{remark}{Remark}

\newtheorem{example}{Example}[section]
\newtheorem{proof}{Proof}

\newtheorem{acknowledgement}{Acknowledgement}

\newcommand{\bt}{\begin{theorem}}
\newcommand{\et}{\end{theorem}}
\newcommand{\bl}{\begin{lemma}}
\newcommand{\el}{\end{lemma}}
\newcommand{\bp}{\begin{proposition}}
\newcommand{\ep}{\end{proposition}}
\newcommand{\bex}{\begin{example}}
\newcommand{\eex}{\end{example}}
\newcommand{\bc}{\begin{corollary}}
\newcommand{\ec}{\end{corollary}}
\newcommand{\bo}{\begin{proof}}
\newcommand{\eo}{\end{proof}}
\newcommand{\bd}{\begin{definition}}
\newcommand{\ed}{\end{definition}}
\newcommand{\br}{\begin{remark}}
\newcommand{\er}{\end{remark}}
\newcommand{\be}{\begin{enumerate}}
\newcommand{\ee}{\end{enumerate}}

\begin{document}

\title{Distal actions and ergodic actions on compact groups}
\author{C. R. E. Raja}
\date{}
\maketitle

\let\epsi=\epsilon 
\let\vepsi=\varepsilon 
\let\lam=\lambda 
\let\Lam=\Lambda 
\let\ap=\alpha 
\let\vp=\varphi 
\let\ra=\rightarrow 
\let\Ra=\Rightarrow 
\let \Llra=\Longleftrightarrow 
\let\Lla=\Longleftarrow 
\let\lra=\longrightarrow 
\let\Lra=\Longrightarrow 
\let\ba=\beta 
\let\ov=\overline 
\let\ga=\gamma 
\let\Ba=\Delta 
\let\Ga=\Gamma 
\let\Da=\Delta 
\let\Oa=\Omega 
\let\Lam=\Lambda 
\let\un=\upsilon

\begin{abstract}
Let $K$ be a compact metrizable group and $\Ga$ be a group of 
automorphisms of $K$.  We first show that each $\ap \in \Ga$ is 
distal on $K$ implies $\Ga$ itself is distal on $K$, a local to global 
correspondence provided $\Ga$ is a generalized $\ov{FC}$-group
or $K$ is a connected finite-dimensional group.  We show that $\Ga$ 
contains an ergodic automorphism when $\Ga$ is nilpotent and ergodic on a 
connected finite-dimensional compact abelian group $K$.  
\end{abstract} 



\begin{section}{Introduction}

We shall be considering actions on compact groups.  By a compact group we 
shall mean a compact metrizable group and by an automorphism we shall 
mean a continuous automorphism.  For a compact group $K$, 
${\rm Aut }~(K)$ denotes the group of automorphisms of $K$.  
An action of a topological group $\Ga$ on a compact metrizable group $K$ 
by automorphisms, is a homomorphism $\phi \colon \Ga \ra {\rm Aut }~(K)$
such that the map $(\ap , x) \mapsto \phi (\ap )(x)$ is a continuous 
map: when only one action is studied or when there is no confusion 
instead of $\phi (\ap )(x)$ we 
write $\ap (x)$ for $\ap \in \Ga$ and $x\in K$.  In such cases, the map 
$\phi$ is said to define the action of $\Ga$ on $K$ and such actions are 
called algebraic actions.    

We shall assume that a topological group $\Ga$ acts on a compact 
metrizable group $K$.  For each $\ap\in \Ga$, $(n , a) \mapsto \ap ^n (a)$ 
defines a ${\mathbb Z}$-action on $K$ and this action on $K$ is called 
${\mathbb Z}_\ap$-action.  Suppose $K_1 \supset K_2$ are closed 
$\Ga$-invariant subgroups of $K$ such that $K_2$ is normal in $K_1$.  By 
an action of $\Ga$ on $K_1/K_2$, we mean the canonical action of $\Ga$ on 
$K_1/K_2$ defined by $\ap (xK_2) = \ap (x) K_2$ for all $x \in K_1$ and 
all $\ap \in \Ga$.

Suppose $\Ga$ acts on the compact groups $K$ and $L$.  We say that $K$ and 
$L$ are $\Ga$-isomorphic if there exists a continuous isomorphism $\Phi 
\colon K \ra L$ such that $\Phi (\ap (x)) = \ap (\Phi (x) )$ for all 
$\ap \in \Ga$ and $x \in K$.  

It is interesting to find properties of group actions that hold if the 
property holds for every ${\mathbb Z}_\ap$-action.  We term any such 
property a local to global correspondence as this 
property holds for the whole group $\Ga$ when it holds locally at every 
point of $\Ga$.  We first state the follwing well-known classical local to 
global correspondence for linear actions on vector spaces proof of which 
may be found in \cite{CR}.

\vskip 0.2in 

\noindent {\bf Burnside Theorem:}  Let $V$ be a finite-dimensional 
vector space over reals and $G$ be a finitely generated subgroup of 
$GL(V)$, the group of linear transformations on $V$.  If each element of 
$G$ has finite order, then $G$ itself is a finite group.

\vskip 0.2in
 
Main aim of the note is to exhibit such local to global correspondences 
for algebraic actions on compact groups.  

\bd
We say that the action of $\Ga$ on $K$ is distal if for any $x\in K 
\setminus (e)$, $e$ is not in the closure of the orbit $\Ga (x) = \{ \ap 
(x) \mid \ap \in \Ga \}$.  In such case, we say that $\Ga$ is distal (on 
$K$).    
\ed

We now introduce a type of action which is obiviously distal.

\bd
We say that the action of $\Ga$ on $K$ is compact (respectively, finite) 
if the group $\phi (\Ga )$ is contained in a compact (respectively, 
finite) subgroup of ${\rm Aut}~(K)$ where $\phi$ is the map defining the 
action of $\Ga$ on $K$.    
\ed

We now see the notion of ergodic action which is orthogonal to distal 
action.  

\bd
Let $K$ be a compact group and $\omega _K$ be the normalized Haar measure 
on $K$.  We say that an (algebraic) action of $\Ga$ on $K$ is ergodic if 
any $\Ga$-invariant Borel set $A$ of $K$ satisfies $\omega _K(A) =0$ or 
$\omega _K(A) =1$.
\ed

\bd
Let $K$ be a compact group and $\ap$ be a continuous automorphism of $K$.  
If the action of ${\mathbb Z}_\ap$ on $K$ is distal (respectively, 
ergodic), then  we say that $\ap$ is a distal (respectively, ergodic) 
automorphism of $K$ or $\ap$ is distal (respectively, ergodic) on $K$.   
\ed

It is easy to see that $\Ga$ is distal implies each $\ap \in \Ga$ 
is distal.  In general each $\ap \in \Ga$ is distal need not imply $\Ga$ 
is distal (Example 1, \cite{Ro}).  For actions on 
connected Lie groups \cite{Ab} and for certains actions on $p$-adic Lie 
groups \cite{Ra1} the local to global correspondence, that is, from each 
$\ap \in \Ga$ being distal on $K$ to the whole group $\Ga$ being distal on 
$K$ holds: distal notion has canonical extension to actions on 
locally compact spaces \cite {Ei}.  Recently \cite{JR} showed under 
certain conditions on $K$ and $\Ga$ that each $\ap \in \Ga$ is distal and 
the whole group $\Ga$ is distal are equivalent to the action being 
equicontinuous (that is, having invariant neighbourhoods).  

We now introduce a class of groups whose action is one of the main studies 
in this article.  

\bd
A locally compact group $G$ is called a generalized $\ov{FC}$-group 
if $G$ has a series $G=G_0\supset G_1\supset \cdots\supset G_n=\{e \}$
of closed normal subgroups such that $G_i/G_{i+1}$ is a compactly 
generated group with relatively compact conjugacy classes for 
$i=0,1,\cdots, n-1$.
\ed

It follows from Theorem 2 of \cite{L} that compactly generated locally 
compact groups of polynomial growth are generalized $\ov{FC}$-groups and 
any polycyclic group is a generalized $\ov{FC}$-group.  

It can be easily seen that the class of generalized $\ov{FC}$-groups is 
stable under forming continuous homomorphic images and closed subgroups.  
If $H$ is a compact normal subgroup of a locally compact group $G$ such 
that $G/H$ is a generalized $\ov{FC}$-group, then it is easy to see that 
$G$ is also a generalized $\ov{FC}$-group.  

In this article we first invesigate the connection between distal 
actions and ergodic actions and apply it to prove a local to global 
correspondence for distal actions if $\Ga$ is a generalized 
$\ov{FC}$-group.  

We next consider finite-dimensional compact groups and prove the local to 
global correspondence for any distal action on compact connected 
finite-dimensional groups.  We also prove a structure theorem for distal 
actions on compact connected finite-dimensional abelian groups along the 
lines of \cite{CG}. 

The study of ergodic actions on compact groups is a key tool in proving 
the afore-stated results.  Using these methods we prove the existence of 
ergodic automorphim in $\Ga$ when $\Ga$ is nilpotent and ergodic on a 
compact connected finite-dimensional abelian group: when $\Ga$ is abelian 
this result is shown in \cite{Be}. 

Having explained our results, it is easy to see that only $\phi (\Ga)$ 
matters and not all of $\Ga$.  So, we may assume that $\Ga$ is a group of 
automorphisms of $K$.  

\end{section}

\bigskip 
\begin{section}{Distal and Ergodic}

In this section we explore the connection between distal actions and 
ergodic actions on compact (metrizable) groups using the dual structure of 
compact groups. 

Let $K$ be a compact group and $\Ga$ be a group acting on $K$.  Let $\hat 
K$ be the equivalent classes of 
continuous irreducible unitary representations of $K$.  If $\pi$ is a 
continuous irreducible unitary representation of $K$, then $[\pi]\in \hat 
K$ denotes the set of all continuous irreducible unitary representations 
of $K$ that are unitarily equivalent to $\pi$.  We write $\pi _1 \sim \pi 
_2$ if $\pi_1, \pi _2 \in [\pi]$ for some $[\pi] \in \hat K$.  For a 
continuous irreducible unitary representation $\pi$ of $K$ and 
$\ap \in \Ga$, $\ap (\pi )$ is defined by 
$$\ap (\pi ) (x) = \pi (\ap ^{-1}(x) )$$ for all $x \in K$ and it can 
be easily verified that $\ap (\pi )$ is also a  continuous irreducible 
unitary representation of $K$.  If $\ap \in \Ga$ and 
$\pi_1, \pi _2 \in [\pi]$, then $\ap (\pi _1) \sim \ap (\pi _2)$.  
Thus, the map $(\ap, ~ [\pi ])\mapsto \ap [\pi ]= [\ap (\pi )]$ is 
well-defined and is known as the dual of action of $\Ga$ on the dual $\hat 
K$ of $K$.  For $k \geq 1$, let $U_k({\mathbb C})$ be the group of 
unitaries on ${\mathbb C}^k$ and $I_k$ denote the identity matrix in 
$U_k({\mathbb C})$.  Then $U_k({\mathbb C})$ is a compact group and for 
each $[\pi ] \in \hat K$, there exists a $k\geq 1$ such that 
$\pi (x) \in U_k({\mathbb C})$ for all $x \in K$: see \cite{Fo} for 
details on representations of compact groups.  

\bp\label{e}
Let $K$ be a (non-trivial) compact group and $\Ga$ be a group of 
automorphisms of $K$.  Then the following are equivalent: 
\be
\item [(1)] $\Ga$ is distal on $K$;

\item [(2)] for each $\Ga$-invariant non-trivial closed subgroup $L$ of 
$K$, action of $\Ga$ on $L$ is not ergodic;

\item [(3)] for each $\Ga$-invariant non-trivial closed subgroup $L$ of 
$K$, there exists a non-trivial continuous irreducible unitary 
representation $\pi$ of $L$ such that the orbit 
$\Ga [\pi ] = \{ \ap [\pi ] \mid \ap \in \Ga \}$ is finite in $\hat L$.
\ee
\ep

\bo
Let $L$ be a non-trivial $\Ga$-invariant closed subgroup of $K$.  If the 
action of $\Ga$ on $L$ is ergodic, then by Theorem 2.1 of \cite{Be}, $\Ga 
(x)= \{ \ap (x) \mid \ap \in \Ga \}$ is dense in $L$ for some $x \in L$.  
Since 
$L$ is non-trivial, $x\not = e$ and hence $e$ is in the closure of
$\Ga (x)$ for $x \not = e$.  Thus, we get that (1) $\Ra$ (2) and 
that (2) $\Ra$ (3) follows from Theroem 2.1 of \cite{Be}.  

Now assume that (3) holds.  Let $x\not = e$ be in $K$ and $L$ be the 
closed subgroup generated by $\Ga (x)$.  Then $L$ is a non-trivial 
$\Ga$-invariant closed subgroup of $K$.  Then by assumption there exists a 
non-trivial $[\pi _1]\in \hat L$ such that $\Ga ([\pi _1])$ is finite.  
Let $\Ga _0 = \{ \ap \in \Ga \mid \ap (\pi _1) \sim \pi _1 \}$.  
Then $\Ga _0$ is a closed subgroup of $\Ga$ of 
finite index.  Let $\Ga _1 = \cap _{\ap \in \Ga }\ap \Ga _0 \ap ^{-1}$.  
Then $\Ga _1$ is a normal subgroup of $\Ga$ of finite index and $\Ga _1$ 
is contained in $\Ga _0$.  Let $A= \{ [\pi] \in \hat L \mid \Ga _1 [\pi ] 
= [\pi ] \}$.  Then $A$ contains $\pi _1$.  
Since $\Ga _1$ is normal in $\Ga$, $A$ is $\Ga$-invariant.  
Let $L_1 = \cap _{[\pi ] \in A} \{ x \in L \mid \pi (x) = \pi (e) \}$.  
Then $L_1$ is a $\Ga$-invariant closed normal subgroup of $L$ and $L_1$ 
is a proper subgroup of $L$ as $A$ is non-trivial.  If $e$ is in the 
closure of $\Ga (x)$, then since $\Ga /\Ga _1$ is finite, $e$ is  in the
closure of $\Ga _1 (x)$.  Let $\ap _n \in \Ga _1$ be such that $\ap _n (x 
) \ra e$ and $[\pi ] \in A$.  Then there exist $u_n \in U_k ({\mathbb C})$ 
($k$ may depend on $\pi$) such that $$ u_n^{-1}\pi (g) u_n  = 
\pi (\ap _n (g))$$ for all $g \in L$.  This implies that 
$$ u_n^{-1}\pi (x) u_n = \pi (\ap _n (x)) \ra \pi (e) = I_k$$ as $n \ra 
\infty$.  Since $U_k ({\mathbb C})$ is compact, $\pi (x)=I_k$.  
This implies that $x \in L_1$ which is a contradiction as 
$L_1$ is a proper $\Ga$-invariant subgroup of $L$ and $L$ 
is the closed subgroup generated by $\Ga (x)$.  Thus, $e$ is not in the 
closure of $\Ga (x)$.  Hence (3) $\Ra $ (1).  
\eo

We now prove a result for compact abelian groups by employing the dual 
structure of locally compact abelian groups.  For a locally compact 
abelian group $G$, a continuous homomorphism of $G$ into the circle 
group $T=\{ z\in {\mathbb C} \mid |z| =1 \}$ is known as a character of 
$G$ and the dual group of $G$ denoted by $\hat G$ is defined to be the 
group of all characters of $G$.  The group $\hat G$ is a locally compact 
abelian group with the standard compact-open topology and the dual of 
$\hat G$ is (isomorphic to ) $G$.  

It is known that $G$ is compact if and only if 
$\hat G$ is discrete and there is a one-one correspondence 
between closed subgroups of $G$ and the quotients of $\hat G$: 
\cite{Mor} for deatils on duality of locally compact abelian groups.  

If $G$ is a group and $A_1, A_2, \cdots , A_n$ are subsets of $G$, then 
$<A_1, \cdots , A_n >$ denotes the group generated by 
$A_1, A_2, \cdots , A_n$ and if any $A_i = \{ g \}$, we may write 
$g$ instead of $\{ g \}$.  

\bl\label{al}
Let $K$ be a compact abelian group and $\Ga$ be a group of automorphisms 
of $K$.  Let $\ap$ be an automorphism of $K$ such that $\ap \Ga \ap ^{-1} 
= \Ga$.  Suppose the action of $\Ga$ is not ergodic on $K$ and for each 
$\ap$-invariant proper closed subgroup $L$ of $K$, the action of ${\mathbb 
Z}_\ap$ on $K/L$ is not ergodic.  Then there exists a non-trivial 
character  $\chi$ on $K$ such that the orbit $\{ \ba (\chi ) \mid \ba \in 
<\Ga, \ap > \}$ is finite - in other words, the group generated by $\Ga$ 
and $\ap$ is not ergodic on $K$.
\el

\bo
We first note that the assumption on $\ap$ is equivalent to saying that 
for any $\ap$-invariant non-trivial subgroup $A$ of $\hat K$ there exists 
a non-trivial character $\chi \in A$ such that the orbit $\{ \ap ^n (\chi 
) \mid n \in {\mathbb Z} \}$ is finite.  

Let $A = \{ \chi \in \hat K \mid \Ga (\chi ) ~~{\rm is ~~ finite} \}$.  
Since $\Ga$ is not ergodic on $K$, $A$ is non-trivial.  Since $\ap \Ga \ap 
^{-1} =\Ga$, $A$ is $\ap$-invariant.  By assumption on $\ap$, there exists 
a non-trivial $\chi _0$ in $A$ such that $\ap ^k (\chi _0 ) = \chi _0$ for 
some $k \geq 1$.  Then $$\Ga \ap ^n (\chi _0 ) \subset \cup _{i =1}^k \Ga 
\ap ^i (\chi _0)$$ for all $n \in {\mathbb Z}$.  
Since $\chi _0\in A$ and $A$ is $\ap$-invariant, we get that  
$\{ \ba (\chi _0) \mid \ba \in <\Ga, \ap > \}$ is finite.
\eo

\bl\label{l0} 
Let $K$ be a compact abelian group and $\Ga$ be a group of 
automorphisms of $K$.  Suppose $\Ga$ is a generalized $\ov{FC}$-group and 
for each $\ap \in \Ga$ and each $\ap$-invariant proper closed 
subgroup $L$ of $K$, the action of ${\mathbb Z}_\ap$ on $K/L$ is not 
ergodic.  Then there exists a non-trivial character $\chi$ on $K$ such 
that the corresponding $\Ga$-orbit $\{\ap (\chi ) \mid \ap \in \Ga \}$ is 
finite or equivalently the action of $\Ga$ on $K$ is not ergodic.
\el

\bo
Since $K$ is compact abelian, ${\rm Aut}~(K)$ is totally disconnected and 
hence by Proposition 2.8 of \cite{JR}, $\Ga$ contains a compact open 
normal 
subgroup $\Da$ such that $\Ga/\Da$ contains a polycyclic subgroup of 
finite index.  Let $\Lam$ be a closed normal subgroup of $\Ga$ of finite 
index containing $\Da$ such that $\Lam /\Da$ is polycylic.  
Let $\Lam _0 = \Lam $ and $\Lam _i= [\Lam _{i-1}, \Lam _{i-1}]$ for 
$i\geq 1$.  Then there exists a $k \geq 0$ such that $\Lam _k\Da \not = 
\Da$ and $\Lam_{k+1} \Da =\Da$.  It can be easily seen that each 
$\Lam _i \Da$ is finitely generated modulo $\Da$.  For $0\leq i\leq k$, 
let $\ap _{i,1}, \cdots , \ap _{i,m} $ be in $\Lam _i$ such that 
$\ap _{i,1}, \cdots , \ap _{i,m} $ and $\Da\Lam _{i+1}$ generate 
$\Da\Lam _i$.  It can be easily seen that $\ap _{i, j}$ 
normalizes $<\ap _{i, 1}, \cdots , \ap _{i, j-1}, \Lam _{i+1}, \Da >$ 
for all $i$ and $j$.  Then repeated applicaton of 
Lemma \ref{al} yields a non-trivial character $\chi \in \hat K$ 
such that the orbit $\Lam (\chi )$ is finite.  Since $\Lam$ is a 
normal subgroup of finite index in $\Ga$, $\Ga (\chi )$ is also finite.  
\eo

We next prove a lemma which shows that the (global) distal condition in 
Proposition \ref{e}, can be relaxed to the local distal condition provided 
$\Ga$ is a generalized $\ov{FC}$-group.  

\bl\label{l1} Let $K$ be a compact abelian group and $\Ga$ be a group of 
automorphisms of $K$.  Suppose $\Ga$ is a generalized $\ov{FC}$-group and 
each $\ap \in \Ga$ is distal on $K$.  Then $\Ga$ is not ergodic on 
$K$ or equivalently there exists a non-trivial character $\chi$ on $K$ 
such that the corresponding $\Ga$-orbit $\{\ap (\chi ) \mid \ap \in \Ga 
\}$ is finite.  
\el

\bo
Let $\ap \in \Ga$ and $L$ be a $\ap$-invariant proper closed subgroup of 
$K$.  Since $\ap$ is distal on $K$, the action ${\mathbb Z}_\ap$ on $K/L$ 
is also distal (Corollary 6.10 of \cite{BJM}).  This shows by Proposition 
\ref{e} that the action of ${\mathbb Z}_\ap$ is not ergodic on $K/L$.  
Thus, the result follows from Lemma \ref{l0}.  
\eo
\end{section}

\begin{section}{Distal actions}

In this section we prove that the distal condition has local to global 
correspodence for actions on compact groups provided the group of 
automorphisms is a generalized $\ov{FC}$-group.

\bt\label{lg}
Let $K$ be a compact group and $\Ga$ be a group of automorphisms of $K$.  
Suppose $\Ga$ is a generalized $\ov{FC}$-group.  Then the following are 
equivalent: 

\be
\item  each $\ap \in \Ga$ is distal on $K$;

\item the action of $\Ga$ on $K$ is distal.  
\ee
\et

\bo
Suppose each $\ap \in \Ga$ is distal on $K$.  Let $x\in K$ be such 
that $e$ is in the closure of the orbit $\Ga (x)$.  We now 
claim that $x= e$. 

\noindent {\bf Case (i):}  Suppose $K$ is abelian.  Let $L$ be a 
non-trivial $\Ga$-invariant closed subgroup of $K$.  It follows from 
Lemma \ref{l1} that $\Ga$ is not ergodic on $L$.  Since $L$ is arbitrary 
$\Ga$-invariant closed subgroup, by Proposition \ref{e} we get that 
$\Ga$ is distal on $K$.  

\noindent {\bf Case (ii):}  Suppose $K$ is connected.  Let $T$ be a 
maximal compact connected abelian subgroup of $K$ containing $x$ (see 
\cite{Hof}).  Then ${\rm Aut}(K) = {\rm Inn}(K) \Oa$ where $\Omega = \{ 
\ap \in {\rm Aut}(K)\mid \ap (T) = T\}$ (see \cite{Hof}).  
Let $\Ga ' = \Ga {\rm Inn}(K)$ and $\Omega '= (\Ga '\cap \Omega )$.  Then 
$\Ga '$ and $\Oa '$ are also generalized $\ov{FC}$-groups.  
Since ${\rm Aut}(K) = {\rm Inn}(K) \Oa$, 
$\Ga ' = {\rm Inn}(K) \Omega '$.  Since $e$ is in the closure of 
$\Ga (x)$, $e$ is in the closure of $\Ga ' (x)$.  Since ${\rm Inn }~(K)$ 
is compact, $e$ is in the closure of $\Oa ' (x)$.  As $x \in T$, applying 
case (i), we get that $x=e$.  

\noindent {\bf General Case:}  Let $K$ be any compact group and $K_0$ 
be the connected component of $e$ in $K$.  Then $K_0$ is $\Ga $-invariant 
and 
by Corollary 6.10 of \cite{BJM}, each $\ap \in \Ga$ is distal on $K/K_0$.  
Since $K/K_0$ is totally disconnected, by Proposition 2.8 and 
Lemma 2.3 of \cite{JR}, $K/K_0$ has arbitrarily small 
compact open subgroups invariant under $\Ga$.  This shows that $x \in 
K_0$.  Now $x=e$ follows from case (ii).  
\eo

As a consequence of results proved so far, we now prove an initial result 
on the existence of ergodic automorphisms in $\Ga$ when the action of 
$\Ga$ on $K$ is ergodic.  

\bp
Let $K$ be a compact group and $\Ga$ be a group of automorphisms of $K$.  
Suppose $\Ga$ is a generalized $\ov{FC}$-group and the action of $\Ga$ on 
$K$ is ergodic.  Then we have the following:

\be
\item [(i)] there exist a $\ba \in \Ga$ and a $\ba$-invariant non-trivial 
closed subgroup $L$ of $K$ such that the action of ${\mathbb Z}_\ba$ on 
$L$ is ergodic.

\item [(ii)] In addition if $K$ is abelian, there exist a $\ap \in \Ga$ 
and a $\ap$-invariant proper closed subgroup $L$ of $K$ such that the 
action of ${\mathbb Z}_\ap$ on $K/L$ is ergodic;
\ee
\ep

\bo
Suppose for each $\ap \in \Ga$ and each $\ap$-invariant non-trivial 
closed subgroup $L$ of $K$, the action of ${\mathbb Z}_\ap$ on $L$ is not 
ergodic.  Then by Proposition \ref{e}, each $\ap \in \Ga$ is distal on 
$K$.  By Theorem \ref{lg}, the action of $\Ga$ on $K$ is distal and 
hence by Proposition \ref{e}, the action of $\Ga$ on $K$ is not ergodic.  
Thus, (i) is proved.   

We now assume that $K$ is abelian.  Suppose for every $\ap \in \Ga$ and 
for every proper closed $\ap$-invariant subgroup $L$ of $K$, the action of 
${\mathbb Z}_\ap$ on $K/L$ is not ergodic.  By Lemma \ref{l0}, the action 
of $\Ga$ on $K$ is not ergodic.  Thus, (ii) is proved.  
\eo

The following example shows that ergodic action of general, even a 
commutative group $\Ga$ on a compact abelian group need not imply 
the existence of a non-trivial subgroup or a non-trivial quotient that  
admits an ergodic ${\mathbb Z}_\ap$-action for some $\ap \in \Ga$.  

\bex
Let $(F_n)$ be an strictly increasing sequence of finite groups 
(one may take $F_n = \prod _{k=1} ^n {\mathbb Z}/k{\mathbb Z}$ ) and $A= 
\cup F_n$.  Let $K = M^A$ where $M$ is a compact abelian group.  
Then $K$ is a compact abelian group 
whose dual $\hat K$ consists of functions $f \colon A \ra \hat M$ 
such that $f(A)$ is finite (see 
Theorem 17 of \cite{Mor}).  We consider the shift action of $A$ on $K$ 
defined by $ag(x) = g(a^{-1}x)$ for all $g\in K= M^A$ and all $a, x\in A$.  
Then the dual 
action of $A$ on the dual $\hat K$ is given by $af(b) = f(a^{-1}b)$ for 
all $f \in \hat K$ and all $a, b\in A$.  Let $f\in \hat K$ be 
non-trivial.  Then there 
exists a $F_n$ such that $f(a)$ is the trivial character on $M$ for all 
$a \not \in F_n$.  For $i\geq 1$, let $a_i \in 
F_{n+i}\setminus F_{n+i-1}$ and $a \in F_n$ be such that $f(a)$ 
is not the trivial character.  Then $a_j^{-1}a_i a \not \in F_n $ if 
$i\not = j$ and hence $a_i^{-1} a_j f(a) = f(a_j^{-1} a_i a)$ is the 
trivial character if $i\not = j$.  This shows that  $a_i f\not =a_j f$ 
for all $i\not = j$.  Thus, the orbit $Af$ is infinite for any 
non-trivial $f \in \hat K$.  This implies that the action of $A$ on $K$ is 
ergodic.  Since any $a \in A$ has finite order, the action of 
${\mathbb Z}_a$ is never ergodic for any $a\in A$.
\eex
\end{section}

\begin{section}{Finite-dimensional compact groups}

In this section we consider finite-dimensional compact groups.  
Let ${\mathbb Q}_d^r$ be the group ${\mathbb Q}^r$ with discrete topology.  
We may regard ${\mathbb Q}_d^r$ as a finite-dimensional vector space over 
${\mathbb Q}$.  Let $B_r$ denote the dual of ${\mathbb Q}_d^r$.  Then 
$B_r$ is a compact connected group of finite-dimension and any compact 
connected finite-dimensional abelian group is a quotient of $B_r$ for some 
$r$: see \cite{Mor}.  

We first show that distal condition for algebraic actions on $B_r$ has 
local to global correspondence with no restriction on the acting group 
$\Ga$.  
The dual of any automorphism of $B_r$ is a ${\mathbb Q}$-linear 
transformation of ${\mathbb Q}_d^r$ onto 
${\mathbb Q}_d^r$.  It can be easily seen that any group of unipotent 
transformations of ${\mathbb Q}_d^r$ is distal on $B_r$.  We now show 
that upto finite extensions these are the only distal actions on $B_r$ 
which is a structure theorem for distal actions on $B_r$ along the lines 
of \cite{CG}.  

\bp\label{fd}
Let $\Ga$ be a group of automorphisms of $B_r$.  Suppose $\Ga$ is 
distal on $B_r$.  Then $B_r$ has a series 
$$B_r= K_0 \supset K_1 \supset 
K_2 \supset \cdots \supset K_{n-1} \supset K_n =(e)$$ of closed connected 
$\Ga$-invariant subgroups such that the action of 
$\Ga$ on $K_i/K_{i+1}$ is finite for any $i\geq 0$.  In particular, $\Ga$ 
is a finite extension of a group of unipotent transformations of 
${\mathbb Q}_d^r$.   
\ep

\bo
Let $\Ga$ act distally on $B_r$.  Then by Proposition \ref{e}, there 
exists a non-trivial $\chi _1\in {\mathbb Q}_d^r$ such that orbit $\Ga 
(\chi _1)$ is finite.  Let $\tilde \Ga _ 1= \{ \ap \in \Ga \mid 
\ap (\chi _1 ) = \chi _1\}$.  Then $\tilde \Ga _1$ is a subgroup of 
finite index in $\Ga$.  Let $\Ga _1$ be a normal subgroup of finite index 
in $\Ga$ and contained in $\tilde \Ga _1$.  Let $A_1 = \{ \chi \in 
{\mathbb Q}_d^r \mid \Ga _1 (\chi ) = \chi \}$.  Then $A_1$ is a 
non-trivial $\Ga$-invariant ${\mathbb Q}$-vector subspace of 
${\mathbb Q}_d^r$.  Let $K_1$ be the closed 
subgroup of $K$ such that the dual of $K/K_1$ is $A_1$.  Then $K_1$ is 
a proper $\Ga$-invariant closed subgroup of $K$ and the action of $\Ga$ on 
$K/K_1$ is finite.  Since the dual of $K_1$ is the ${\mathbb Q}$-vector 
space ${\mathbb Q}_d^r/A_1$, $K_1\simeq B_{r_1}$ for some $r_1 <r$.  
If $K_1 \not = (e)$, get $K_2$ by applying the above process to 
$K_1$ in place of $K$.  Since $B_r$ has finite-dimension and each 
$K_i$ is connected, proceeding this way we obtain a series $$K_0= B_r 
\supset K_1 \supset K_2 \supset \cdots \supset K_{n-1} \supset K_n =(e)$$ 
of $\Ga$-invaraint closed connected subgroups such that 
the action of $\Ga$ on $K_i/K_{i+1}$ is finite.  
\eo

\bt\label{lg3}
Let $\Ga$ be a group of automorphisms of $B_r$.  Suppose each 
$\ap \in \Ga$ is distal on $B_r$.  If the dual action of $\Ga$ on 
${\mathbb Q}_d^r$ is irreducible, then $\Ga$ is finite.  In general, $\Ga$ 
is distal on $B_r$.  
\et

\bo
Let $\ap \in \Ga$.  By considering the dual action of $\ap$, we may view 
$\ap$ as a linear map on ${\mathbb Q}_d^r$.  Then by Proposition \ref{fd}, 
eigenvalues of $\ap$ are of absolute value one.  If the dual action of 
$\Ga$ on ${\mathbb Q}_d^r$ is irreducible, then let $V = {\mathbb 
Q}_d^r \otimes {\mathbb R}$ be the corresponding 
vector space over ${\mathbb R}$.  Then $V$ is also $\Ga$-irreducible. 
By \cite{CG}, $\Ga$ is contained in a compact subgroup of 
$GL(V)$.  By Proposition \ref{fd}, $\ap ^k$ is unipotent for 
some $k \geq 1$.  Since $\ap$ is in a compact group, we get that $\ap$ 
is of finite order.  Thus, every element $\ap$ of $\Ga$ has finite order.  
It follows from Lemma 4.3 of \cite{Be} that $\Ga$ is finite.  
\eo

We now proceed to show that ergodic action of $\Ga$ on a 
finite-dimensional compact connected abelian group yields an ergodic 
automorphism in $\Ga$ provided $\Ga$ is nilpotent. 

\bl\label{d}
Let $\Ga$ be a group of automorphisms of a compact group $K$ and $\ap$ 
be an automorphism of $K$.  Suppose $\Ga$ and $\ap$ are distal on $K$ and 
$\ap \Ga \ap ^{-1} =\Ga$.  Then the group generated by $\Ga$ and $\ap$ is 
distal on $K$. 
\el

\bo
Let $\Da$ be the group generated by $\Ga$ and $\ap$.  Let $L$ be a closed 
subgroup of $K$ invariant under $\Da$.  Let $A= \{ [\pi ] \in \hat L \mid 
\Ga ([\pi ])~~{\rm is~~finite} \}$.  Since $\Ga$ is normalized 
by $\ap$, $A$ is $\ap$-invariant.  Since $\ap$ and $\Ga$ are distal on 
$K$, it follows from Proposition \ref{e} that there exists a non-trivial 
$[\pi _0]\in A$ such that 
$\ap ^k (\pi _0)  \sim \pi _0$ for some $k \geq 1$.  Now, 
$\Ga \ap ^i [\pi _0 ] \subset \cup _{j=1}^k \Ga (\ap ^j [\pi _0])$ for 
any $i \in {\mathbb Z}$.  This implies that the orbit $\Da [\pi _0]$ is 
finite.  This shows by Proposition \ref{e} that $\Da$ is distal on $K$.
\eo

\bl\label{L3}
Let $\ap$ be an ergodic automorphism of $B_r$ and $L$ be a closed 
connected subgroup of $B_r$.  Then $\ap$ is ergodic on $L$. 
\el

\bo
It can be easily seen that $\ap$ is ergodic on $B_r$ if and only if no 
root of unity is an eigenvalue of $\ap $ on ${\mathbb Q}_d^r$.  Let $V$ be 
the ${\mathbb Q}$-subspace of ${\mathbb Q}_d^r$ such that the 
dual of $L$ is ${\mathbb Q}_d^r /V$.  Since $\ap$ is ergodic, no root of 
unity is an eigenvalue for $\ap$ on ${\mathbb Q}_d^r$ and hence no root of 
unity is an eigenvalue for $\ap$ on ${\mathbb Q}_d^r/V$.  Thus, $\ap$ is 
ergodic on $L$. 
\eo 

\bl\label{L2}
Let $\ap$ and $\ba$ be automorphisms on $B_r$.  Suppose $\ap$ is contained 
in a group $\Ga$ of automorphisms of $B_r$ such that $\Ga$ is distal and 
$\ba$ is ergodic and normalizes $\Ga$.  Then $\ap^i \ba ^j$ and $\ba ^j 
\ap ^i$ are ergodic for all $i$ and $j$ in ${\mathbb Z}$ with $j \not = 
0$.  
\el

\bo
It is enough to show that $\ap \ba$ and $\ba \ap$ are ergodic.  We first 
prove the case when $\Ga$ is finite.  Assume $\Ga$ is finite.  Let $\chi$ 
be a character such that the orbit $\{ (\ap \ba )^n (\chi ) \mid n \in 
{\mathbb Z} \}$ is finite.  Since $\Ga$ is finite and $\Ga $ is normalized 
by $\ap \ba$, the orbit $\tilde \Ga (\chi )$ is also finite where 
$\tilde \Ga $ is the group generated by $\ap \ba$ and $\Ga$.  Since $\ba 
\in \tilde \Ga$ and $\ba$ is ergodic, we get that $\chi$ is trivial.  
Thus, $\ap \ba$ is ergodic.

We now consider the general case.  Let $V= \{ \chi \in {\mathbb Q}_d^r 
\mid \Ga (\chi )~~{\rm is ~~finite} \}$.  Since $\Ga$ is distal, $V$ is 
a nontrivial ${\mathbb Q}$-subspace and $V$ is invariant under $\ba$ as 
$\Ga$ is normalized by $\ba$.  Let $L$ be the closed connected subgroup 
of $B_r$ such that the dual $B_r/L$ is $V$.  Then $L$ is a proper 
closed connected subgroup invariant under $\Ga$ and $\ba$ and $\Ga$ is 
finite on $B_r/L$.  Then $\ap \ba$ is ergodic on $B_r/L$.  Since the dual 
of $L$ is $ {\mathbb Q}_d^r /V$, $L \simeq B_s$ for some $s <r$.  By Lemma 
\ref{L3}, $\ba$ is ergodic on $L$ and hence by induction on dimension of 
$B_r$, $\ap \ba$ is ergodic on $L$.  Thus, $\ap \ba$ is ergodic on $B_r$.  
Similarly we may show that $\ba \ap$ is also ergodic on $B_r$.  
\eo

\bl\label{L4}
Let $\ap$ be an automorphism of $B_r$.  Then there exists a compact 
connected subgroup $K$ of $B_r$ isomorphic to $B_s$ for some $s>0$ such 
that $\ap$ is ergodic on $K$ and 
$\ap$ is distal on $B_r /K$.  Moreover, if $\Ga$ is a nilpotent group of 
automorphisms of $B_r$ containing $\ap$, then $K$ is $\Ga$-invariant. 
\el

\bo
Let $V_1$ be the ${\mathbb Q}$-subspace of ${\mathbb Q}_d^r$ defined by 
$$V_1 = \{ \chi \in {\mathbb Q}_d^r \mid (\ap ^n(\chi))~~{\rm is ~~ 
finite} \}$$ and define $V_i$ inductively by 
$$V_i = \{ \chi \in {\mathbb Q}_d^r \mid (\ap ^n(\chi) +V_{i-1})~~{\rm is 
~~ finite ~~ in }~~ {\mathbb Q}_d^r /V_{i-1} \}$$ for any $i> 1$.  Then 
each $V_i$ is a $\ap $-invariant ${\mathbb Q}$-subspace.  
Since ${\mathbb Q}_d^r$ has 
finite-dimension over ${\mathbb Q}$, there exists a $n$ such that 
$V_n = V_{n+i}$ for all $i\geq 0$ and for any non-trivial $\chi \in 
{\mathbb Q}_d^r/V_n$,  the orbit $(\ap ^n(\chi) +V_{n})$ is infinite.  Let 
$K$ be a closed subgroup of $B_r$ such that the dual of $K$ is 
${\mathbb Q}_d^r/V_n$.  Then $K$ is $\ap$-invariant and connected.  
The choice of $V_n$ shows that $\ap$ is ergodic on $K$ and $\ap$ is 
distal on $B_r/K$.  

Suppose $\Ga$ is a nilpotent group containing $\ap$.  We show that $V_1$ 
is $\Ga$-invariant by induction on the length of the series $\Ga = \Ga _0 
\supset \cdots \supset \Ga _k=[\Ga , \Ga _{k-1}]\supset \Ga _{k+1}=(e)$.  
Now for $\ba \in \Ga _k$ and $i\in {\mathbb Z}$, 
$\ap ^i \ba = \ba \ap ^i$ and hence $V_1$ is $\Ga _k$-invariant.  
If $V_1$ is $\Ga _{k-j}$-invariant, then for 
$\ba \in \Ga _{k-j-1}$ and $i\in {\mathbb Z}$, $\ap ^i \ba = 
\ba \ap ^i \ba _i$ for some $\ba _i \in \Ga _{k-j}$.  Since $V_1$ is 
$\Ga _{k-j}$-invariant and an iterate of $\ap $ is trivial on $V_1$, we 
get that $V_1$ is $\Ga _{k-j-1}$-invariant.  Hence by induction on $k$, we 
get that $V_1$ is $\Ga$-invariant.  Since each $V_i/V_{i-1}$ is the space 
of all characters whose orbit is finite in ${\mathbb Q}_d^r/V_{i-1}$, we 
get that $V_i$ is $\Ga$-invariant for any $i\geq 1$.
\eo

\bl\label{L5}
Let $\Ga$ be a nilpotent group of automorphisms of $B_r$ and $\ap , \ba 
\in \Ga$.  Let $\Ga _0 = \Ga $ and $\Ga _i = [\Ga , \Ga _{i-1}]$ for 
$i\geq 1$.  Let $k\geq 1$ be such that $\ap \in \Ga _{k-1}\setminus 
\Ga _k$.  Suppose $\ap$ is ergodic on $B_r$ and $\Ga _k$ is distal on 
$B_r$.  Then there exists a $i\geq 0$ such that $\ap ^i\ba$ is ergodic on 
$B_r$. 
\el

\bo
We prove the result by induction on the dimension of $B_r$.  If $r =1$, 
then we have nothing to prove.  So, we may assume that $r>1$.  
If $\ap \ba$ is ergodic on $B_r$, then 
we are done.  So, we may assume by Lemma \ref{L4} that there exists a 
closed connected $\Ga$-invariant proper subgroup $K$ of $B_r$ such that 
$\ap \ba$ is ergodic on $K$ and $\ap \ba $ is distal on $B_r/K$.  
Let $\Da$ be the group generated by $\ap\ba$ and $\Ga _k$.  Then by 
Lemma \ref{d}, $\Da$ is distal on $B_r/K$.  Since $\ap$ and $\ba$ commute 
modulo $\Ga _k$, we 
get that $\ap$ normalizes $\Da$.  By Lemma \ref{L2}, $\ap ^i\ba$ is 
ergodic on $B_r /K$ for all $i \geq 2$.  Since $\ap$ is ergodic on 
$K$, induction hypothesis applied to 
$K$ in place of $B_r$ and $\ap ^2 \ba$ in place of $\ba$, we get that 
$\ap ^j \ba$ is ergodic on $K$ for some $j \geq 2$.  Thus, $\ap ^j \ba $ 
is ergodic on $B_r$ for some $j \geq 2$. 
\eo

\bl\label{L6}
Let $\Ga$ be a nilpotent group of automorphisms of $B_r$ and $\ap , \ba 
\in \Ga$.  Let $\Ga _0 = \Ga $ and $\Ga _i = [\Ga , \Ga _{i-1}]$ for 
$i\geq 1$ and $k\geq 1$ be such that $\ap \in \Ga _{k-1}\setminus \Ga _k$.
Let $K$ be a closed $\Ga$-invariant subgroup of $B_r$ 
isomorphic to $B_s$ for some $s \geq 0$ such that 
$\ap$ is ergodic on $K$ and $\ap$ is distal on $B_r /K$.  
If $\ba$ is ergodic on $B_r /K$ and $\Ga _k$ is distal on $B_r$,   
then there exists $j\geq 0$ such that $\ap ^j \ba$ is ergodic on $B_r$.
\el

\bo
Since $\ap$ normalizes $\Ga _k$, by Lemma \ref{d}, the group generated by 
$\ap$ and $\Ga _k$ is distal on $B_r /K$.  Since $\ba$ centralizes $\ap$ 
modulo $\Ga _k$, it follows from Lemma \ref{L2} that $\ap ^i\ba $ is 
ergodic on $B_r/K$ for all $i\geq 0$.  By Lemma \ref{L5}, $\ap ^j \ba$ is 
ergodic on $K$ for some $j \geq 0$.  This shows that for some $j \geq 0$, 
$\ap ^j \ba$ is ergodic on $B_r$.
\eo

\bl\label{L1}
Let $\Ga$ be a nilpotent group of automorphisms of $B_r$.  
Let $\Ga _0 = \Ga $ and $\Ga _i = [\Ga , \Ga _{i-1}]$ for $i\geq 1$.  
Suppose that the action of $\Ga$ on $B_r$ is ergodic.  Then there exist 
a series 
$$(e) = K_0 \subset K_1 \subset K_2 \subset \cdots \subset K_{m-1} 
\subset K_m =B_r$$ of closed connected $\Ga$-invariant subgroups 
with each $K_i \simeq B_{r_i}$ for some $r _i \geq 0$ and automorphisms 
$\ap _1, \ap _2 , \cdots , \ap _m$ in $\Ga$ 
with the following properties for each $i=1, 2, \cdots , m$: 

\be
\item [(i)] if $k_i$ is the smallest integer $k$ for which $\ap _i \not 
\in \Ga _{k}$, then the action of $\Ga _{k_i}$ on $B_r/K_{i-1}$ is distal;

\item [(ii)] the action of ${\mathbb Z}_{\ap _i}$ on $K_i/K_{i-1}$ is 
ergodic;

\item [(iii)] the action of ${\mathbb Z}_{\ap _i}$ on $B_r/K_{i}$ is
distal:
\ee
\el

\bo
For each $\ap \in \Ga$, if the action of ${\mathbb Z}_\ap$ is distal on 
$B_r$, then by Theorem \ref{lg}, the action of $\Ga$ is distal.  This is a 
contradiction to the ergodicity of $\Ga$ by Proposition \ref{e}.  Thus, 
the action of ${\mathbb Z}_\ap$ is not distal for some $\ap \in \Ga$.  

Since $\Ga$ is nilpotent, there exists a $k$ such that $\Ga _k \not = (e)$ 
and $\Ga _{k+1} =(e)$.  Now, choose $\ap _1 \in \Ga _{k_{1}-1} \setminus 
\Ga _{k_1}$ such that the action of ${\mathbb Z}_{\ap _1}$ is not distal 
on $B_r$ but  the action of $\Ga _{k_1}$ is distal on $B_r$.  By Lemma 
\ref{L4}, there exists a non-trivial $\Ga$-invariant closed connected 
subgroup $K_1$ of $B_r$ isomorphic to $B_{r_1}$ for some $r_1 >0$ such 
that $\ap _1$ is ergodic on $K_1$ and the 
action of ${\mathbb Z}_{\ap _1}$ is distal on $B_r/K_1$.  

Let $L=B_r/K_1$.  Then the action of $\Ga$ on $L$ is ergodic and $L \simeq 
B_{s_1}$ for $s_1 <r$ as $K_1$ is non-trivial.  
By applying induction on the dimension of $B_r$, we get 
$\Ga$-invariant closed connected subgroups $(e) = K_0 \subset K_1 
\subset K_2 \subset
\cdots \subset K_{m-1} \subset K_m =B_r$ and automorphisms $\ap _2 , 
\cdots , \ap _m$ satisfying (i) - (iii) for $2\leq i \leq n$.  
\eo

We now consider connected finite-dimensional 
compact abelian groups.  Let $K$ be a connected finite-dimensional 
compact abelian group.  Then 
$${\mathbb Z}^r \subset \hat K \subset {\mathbb Q}_d^r$$ 
and $K$ is a quotient of $B_r$ for some $r\geq 1$.  
Let $\ap$ be an automorphism of $K$.  Then $\ap$ is an automorphism of 
$\hat K$.  Since ${\mathbb Z}^r \subset \hat K$, $\ap$ has a 
canonical extension to an 
invertible ${\mathbb Q}$-linear map on ${\mathbb Q}_d^r$, 
say $\tilde \ap$.  Thus, any automorphism $\ap$ of $K$ can be lifted to a 
unique automorphism $\tilde \ap$ of $B_r$.  Let $\Ga$ be a group of 
automorphisms of $K$ and $\tilde \Ga$ be the group consisting of lifts 
$\tilde \ap$ of automorphisms $\ap \in \Ga$.  We consider $\Ga$ and 
$\tilde \Ga$ as topological groups with their respective compact-open 
topologies as automorphism groups of $K$ and $B_r$.  By looking at the 
dual action, we can see that the topological groups $\Ga$ and $\tilde \Ga$ 
are isomorphic.  If $\phi \colon B_r \ra K$ is the canonical quotient map, 
then for $\ap \in \Ga$ and $x \in B_r$, we have 
$$\phi (\tilde \ap (x) ) = \ap ( \phi (x) )$$ where $\tilde \ap$ is the 
lift of $\ap $ on $B_r$.   

\bp\label{P1}
Let $K$ be a connected finite-dimensional 
compact abelian group.  Let $\Ga$ be a 
group of automorphisms of $K$ and $\tilde \Ga$ be the corresponding group 
of automorphisms of $B_r$.  Then $\Ga$ is distal (respectively, ergodic) on 
$K$ if and only if $\tilde \Ga$ is distal (respectively, ergodic) on 
$B_r$.
\ep

\bo
For $\chi \in {\mathbb Q}_d^r$, there exists $n \geq 1$ such that $n\chi 
\in \hat K$ and since $\hat K\subset {\mathbb Q}_d^r$, $\Ga$ is ergodic on  
$K$ if and only if $\tilde \Ga$ is ergodic on $B_r$.  Since $K$ is a 
quotient of $B_r$, $\tilde \Ga$ is distal on $B_r$ implies $\Ga$ is distal 
on $K$ (see \cite {BJM}, Corrolary 6.10).  

Suppose $\Ga$ is distal on $K$.  
By Proposition \ref{e}, there exists a non-trivial character 
$\chi _1$ in $\hat K\subset {\mathbb Q}_d^r$ such that $\Ga (\chi _1)$ is 
finite.  Let $V_1 = \{ \chi \in {\mathbb Q}_d^r \mid \tilde \Ga (\chi ) 
~~{\rm is ~~ finite}\}$.  Then $V_1$ is a non-trivial $\tilde 
\Ga$-invariant ${\mathbb Q}$-subspace as $\chi _1 \in V_1$.  Let $M$ be 
a closed subgroup of $B_r$ such that the dual of $M$ is 
${\mathbb Q}_d^r/V_1$.  Then $M$ is $\tilde \Ga$-invariant and $M \simeq 
B_s$ for $s<r$.   
Let $A=V_1\cap \hat K$ and $L$ be a closed subgroup of $K$ such that the 
dual of $L$ is $\hat K / A$.  Then $L$ is $\Ga$-invariant.  Since 
$\hat K/A \subset {\mathbb Q}_d^r/V_1$, $\hat K /A$ has no element of 
finite order and hence $L$ is connected (see Theorem 30 of \cite{Mor}).  
It can be verified that the dimension of $L$ is same as the dimension of 
$M$.  
Hence by induction on the dimension of $K$ we get that $\tilde \Ga$ is 
distal on $M$.  Since the action of $\tilde \Ga$ on $B_r/M$ is finite, 
$\tilde \Ga$ is distal on $B_r$. 
\eo

\bt\label{eag}
Let $K$ be a compact connected finite-dimensional abelian group and $\Ga$ 
be a nilpotent group of automorphisms of $K$.  Suppose $\Ga$ is ergodic 
on $K$.  Then there exists a $\ap \in \Ga$ such that $\ap$ is ergodic on $K$.
\et

\br
This result is another local to global correspondence (that is, no $\ap 
\in \Ga$ is ergodic implies $\Ga$ itself is not ergodic).  
\er

\bo
We first assume that $K=B_r$ for some $r\geq 1$.  By Lemma \ref{L1}, 
there are $\Ga$-invariant closed connected subgroups 
$(e) = K_0 \subset K_1 \subset K_2 \subset \cdots \subset K_{m-1} \subset 
K_m =B_r$ with each $K_i \simeq B_{r_i}$ for some $r_i\geq 0$ and 
automorphisms $\ap _1 ,\ap _2 , \cdots , \ap _m$ in $\Ga$ satisfying 
(i) - (iii) of Lemma \ref{L1}.  We may assume that $K_i \not = K_{i-1}$ 
for $1\leq i \leq m$.  We now prove the result by induction on 
$r$.  If $r=1$, we are done.  Induction hypothesis applied to the action of 
$\Ga$ on $B_r /K_1 \simeq B_{r-r_1}$ 
yields $\ba \in \Ga$ such that $\ba$ is ergodic on 
$B_r/K_1$.  By Lemma \ref{L6}, there exists $\ap \in \Ga$ such that $\ap$ 
is ergodic on $B_r$. 

Now let $K$ be a connected finite-dimensional compact abelian group and 
$r$ be the dimension of $K$.  Then $K$ is a quotient of $B_r$.  Let 
$\tilde \Ga$ be the group of lifts of automorphisms of $\Ga$.  Then by 
Proposition \ref{P1}, $\tilde \Ga$ is ergodic on $B_r$.  It follows from 
the previous case that there exists $\ap \in \Ga$ such that the lift 
$\tilde \ap$ of $\ap$ is ergodic on $B_r$.  Another application of 
Proposition \ref{P1} shows that $\ap$ itself is ergodic on $K$.
\eo

We now show that distal condition for algebraic actions on connected
finite-dimensional compact groups has local to global correspondence
with no restriction on the acting group $\Ga$.

\bt\label{lg2}
Let $\Ga$ be a group of automorphisms of a compact connected 
finite-dimensional group $K$.  Suppose each $\ap \in \Ga$ is distal on 
$K$.  Then the action of $\Ga$ on $K$ is distal. 
\et

\bo
If $K$ is abelian, then the result follows from Proposition \ref{P1} and 
Theoroem \ref{lg3}.  
Suppose $K$ is any finite-dimensional compact connected group.  
Let $x \in K$ and $(\ap _n)$ be a sequence in $\Ga$.  Suppose 
$\ap _n (x) \ra e$.  

Let $T$ be a maximal compact connected subgroup of $K$ containing $x$.  
Since $K$ is a connected group, ${\rm Aut}(K) = {\rm Inn}(K) \Oa$ 
where $\Omega = \{ \ap \in {\rm Aut}(K)\mid \ap (T) = T\}$ 
(see \cite{Hof}).  Let $\Lam = {\rm Inn}(K) \Ga$.  Since ${\rm Inn 
}(K)$ is a compact normal subgroup, each $\ap \in \Lam$ is distal on $K$.  
Let $\ap _n = a_n 
\ba _n$ where $a _n \in {\rm Inn}(K)$ and $\ba _n\in \Oa \cap \Lam$ for 
all $n \geq 1$.  Since ${\rm Inn} (K)$ is compact, by passing to a 
subsequence, if necessary, 
we may assume that $\ba _n (x) \ra e$.  Since $T$ is closed in $K$ which 
is of finite-dimension, $T$ is also of finite-dimension (\cite{Na}) and so 
$\Oa \cap \Lam$ is distal on $T$ and hence $x =e$ as $x \in T$ and 
$\ba _n \in \Oa \cap \Lam$.  Thus, the action of $\Ga$ is distal on $K$. 
\eo

We now prove a structure theorem for distal actions on compact connected 
finite-dimensional abelian groups along the lines of \cite{CG}.  

\bp\label{fd1}
Let $K$ be a compact connected finite-dimensional abelian group and $\Ga$ 
be a group of automorphisms of $K$.  Suppose each $\ap \in \Ga$ is distal 
on $K$ or equivalently the action of $\Ga$ on $K$ is distal.  Then 
$K$ has a finite sequence $$K= K_0 \supset K_1 \supset K_2 \supset \cdots 
\supset K_{n-1} \supset K_n =(e)$$ of closed connected $\Ga$-invariant 
subgroups such that the action of $\Ga$ on $K_i/K_{i+1}$ is finite.    
\ep

\bo
Let $K$ be a compact connected finite-dimensional abelian group.  
Let $r$ be the dimension of $K$.  Then $K$ is a quotient of $B_r$.  
Suppose $\Ga$ is a group of automorphisms of $K$ such that the action of 
$\Ga$ on $K$ is distal.  Let $\tilde \Ga$ be the group of automorphisms of 
$B_r$ consisting of lifts of automorphisms in $\Ga$.  By Proposition 
\ref{P1}, $\tilde \Ga$ is distal on $B_r$.  

By Proposition \ref{fd}, $B_r$ has a series $$B_r= K_0 \supset K_1 \supset 
K_2 \supset \cdots \supset K_{n-1} \supset K_n =(e)$$ of 
$\tilde \Ga$-invariant subgroups such that the action of 
$\tilde \Ga$ on $K_i/K_{i+1}$ is finite for $i\geq 0$.  
Let $\phi \colon B_r \ra K$ be 
the canonical projection.  Let $L_i = \phi (K_i)$ for $1\leq i \leq n$.  
Then each $L_i$ is a closed connected $\Ga$-invariant subgroup of $K$ and 
$L_i \supset L_{i+1}$ for $i\geq 0$.  Since lifting of an automorphism of 
$K$ to an automorphism of $B_r$ is unique, the action of $\Ga $ on 
$L_i/L_{i+1}$ is finite for $i\geq 0$.
\eo

The above can be extended to distal actions on compact Lie groups which 
is already proved in \cite{Ab} for finite-dimensional torus.

\bc\label{c1}
Let $K$ be a compact real Lie group and $\Ga$ be a group of automorphisms 
of $K$.  Suppose each $\ap \in \Ga$ is distal on $K$.  Then there exist 
a series $$K =K_1 \supset K_2 \supset \cdots \supset K_{n-1}\supset 
K_{n}=(e)$$ of $\Ga$-invariant closed 
connected normal subgroups of $K$ such that 
the action of $\Ga $ on $K_i/K_{i+1}$ is compact for any $i\geq 0$.  
\ec

\bo
Since $K$ has only finitely many connected components, we may assume that 
$K$ is connected.  Since $K$ is a compact connected Lie group, 
$K = TS$ where $T$ is the connected component of identity in the center of 
$K$ and $S= [K,K]$ is a compact connected semisimple Lie group 
(see \cite{Hoc}).  Thus, $T$ and $S$ are invariant under $\Ga$.  Since $T$ 
is abelian, by Proposition \ref{fd1}, 
there exists a series 
$$T= T_0\supset T_1\supset \cdots \supset T_{n-1}\supset T_n=(e)$$
of closed connected $\Ga$-invariant subgroups 
such that the $\Ga$-action on $T_i/T_{i+1}$ is finite for all $i$.  For 
$0\leq i \leq n$, let $K_i = T_iS$ and $K_{n+1}=(e)$.  Then each $K_i$ is 
a $\Ga$-invariant closed 
connected normal subgroup of $K$ such that $K_i\supset 
K_{i+1}$ for $0\leq i\leq n$.  For $0\leq i\leq n-1$, $K_i/K_{i+1}$ is 
$\Ga$-isomorphic to $T_i/T_i\cap T_{i+1}S$ and hence the $\Ga$-action on 
$K_i/K_{i+1}$ is finite for $0\leq i \leq n-1$.  Since $S$ is a compact 
connected semisimple Lie group, ${\rm Aut}~(S)$ is compact and hence the 
$\Ga$-action on $K_n/K_{n+1}= S$ is compact.    
\eo

We now provide an example to show that nilpotency assumption on the acting 
group $\Ga$ in Theorem \ref{eag} can not be relaxed: it may be noted that 
Theorem \ref{eag} is true with no restriction on the acting group $\Ga$ if 
the compact group $K$ is a two-dimensional torus.  

\bex
Let $\Ga$ be a subgroup of ${\rm GL} (n , {\mathbb Q})$.  
Let $\Ga ^+$ be the semi-direct product of $\Ga $ and ${\mathbb Q}_d^n$ 
with the canonical action of $\Ga$ on ${\mathbb Q}_d^n$.  We define an 
action of $\Ga ^+$ on ${\mathbb Q}_d^{n+1}$ by 
$$(\ap , w) (q_1, \cdots ,q_n , q_{n+1}) = \ap (q_1, \cdots , q_n ) + 
wq_{n+1} +(0, \cdots , 0, q_{n+1})$$ for all $(\ap , w) \in \Ga ^+$ and 
$(q_1, \cdots ,q_n , q_{n+1}) \in {\mathbb Q}_d^{n+1}$: ${\mathbb Q}_d ^n$ 
is identified as a subset of ${\mathbb Q}_d^{n+1}$ via the canonical 
map $(q_1, \cdots , q_n ) \mapsto (q_1, \cdots , q_n , 0)$.  Considering 
the dual action, we get that $\Ga ^+ \subset {\rm Aut} (B_{n+1})$.  
For $z \in {\mathbb Q}_d^n $, $\Ga ^+ (z) = \Ga (z)$ and for 
$z \in {\mathbb Q}_d ^{n+1} \setminus {\mathbb Q}_d^n$, 
$\Ga ^+ (z)$ can be easily seen to be infinite.  Thus, $\Ga $ is ergodic 
on $B_n$ if and only if $\Ga ^+$ is ergodic on $B_{n+1}$.  
For any $\Ga \subset {\rm GL}(n , {\mathbb Q})$, no $\ap \in \Ga ^+$ is 
ergodic on $B_{n+1}$.  For $n \geq 1$, take $\Ga$ to be the group 
generated by $\ap \in {\rm GL} (n , {\mathbb Q})$ that is ergodic on 
$B_n$.  Then $\Ga ^+$ is a solvable group and is ergodic on $B_{n+1}$ 
but no automorphism in $\Ga ^+$ is ergodic on $B_{n+1}$.  
\eex

\end{section}

\begin{section}{Compact abelian groups}

In this section we obtain a general version of Proposition \ref{fd1} for 
compact abelian groups and also provide an example to show that the 
existence of finite sequence in Proposition \ref{fd1} need not be 
true for connected infinite-dimensional compact abelian groups.  

\bp\label{p1}
Let $K$ be a compact abelian group and $\Ga$ be a group of automorphisms 
of $K$.  Suppose $\Ga$ is distal.  Then there 
exists a collection $(K_i)$ of $\Ga$-invariant closed subgroups of $K$ 
such that 
\be
\item [(1)] $K_0 = K$;

\item [(2)] for $i \geq 0$ either $K_{i+1} = (e)$ or $K_{i+1}$ is a proper 
subgroup of $K_i$;

\item [(3)] the action of $\Ga $ on $K_i /K_{i +1}$ is finite 
for any $i \geq 0$;

\ee
\ep

\bo
Let $K_0= K$.  Then by Proposition \ref{e} there exists a non-trivial 
character $\chi _0$ on $K_0$ such that $\{ \ap (\chi _0) \mid \ap \in \Ga 
\}$ is finite.  Let $\tilde \Ga _0 = \{\ap \in \Ga \mid \ap (\chi _0)=
\chi _0 \}$.  
Then $\tilde \Ga _0$ is a subgroup of finite index in $\Ga$.  Let $\Ga _0$ 
be a normal subgroup of $\Ga$ of finite index and $\Ga _0 \subset \tilde 
\Ga _0$.

Let $A_1 = \{ \chi \in \hat K _0 \mid \Ga _0(\chi ) = \chi \}$.  Since 
$\Ga _0$ is normal in $\Ga$, $A_1$ is a $\Ga $-invaraint 
non-trivial subgroup of $\hat K_0$ and hence there exists a 
proper closed subgroup $K_1$ of $K_0$ such that the dual of $K_0/K_1$ is 
$A_1$.  Then $K_1$ is $\Ga$-invariant and the action of $\Ga $ on 
$K/K_{1}$ is finite.   

If $K_1 = (e) $, then take $K_n =(e)$ for all $n \geq 1$.  
If $K_1\not = (e)$, then get $K_2$ by applying the above arguements to 
$K_1$.  Proceeding this way we obtain a collection $(K_i)$ of 
$\Ga$-invariant closed subgroups of $K$ satisfying conditions (1)-(3). 
\eo       

In contrast to the finite-dimensional case we now show by an example that 
the sequence $(K_i)$ in Proposition \ref{p1} need not be finite.  
Let $T_k$ be the $k$-dimensional torus, a product of $k$ copies of the 
circle group.  Let $\ap _k$ be the automorphism of $T_k$ defined by 
$$\ap _k (x_1, x_2, \cdots , x_k) =(x_1x_2\cdots x_k, x_2x_3\cdots x_k, 
\cdots , x_{k-1}x_k , x_k)$$ for all $(x_1, x_2, \cdots ,x_k) \in T_k$.  
For $0\leq j \leq k$, let $M_{k, j} = \{ (x_1, x_2, \cdots , x_{k-j}, e, 
\cdots ,e) \in T_k \}$.  Then each $M_{k,j}$ is $\ap _k$-invariant and 
$\ap _k$ is trivial on $M_{k,j}/M_{k,j+1}$ for $j \geq 0$.    

We first prove the following fact about $T_k$ and $\ap _k$. 

\bl\label{tl}
Let $T_k$ and $\ap _k$ be as above.  Suppose there exists a series 
$$T_k= M_0 \supset M_1 \supset \cdots \supset M_{n-1}\supset M_n =(e)$$ of 
$\ap_k$-invariant closed subgroups such that for $i\geq 0$, the action of 
${\mathbb Z}_{\ap _k}$ on $M_{i}/M_{i+1}$ is finite and $M_{i}/M_{i+1}$ 
is not finite.  Then $n = k$.  
\el

\bo
Let $V$ be the Lie algebra of $T_k$. We first show that $M_{n-1}$ is 
one-dimensional. 
For $0\leq i <n$, let $V_i$ be the Lie subalgebra of $V$ corresponding 
to the Lie subgroup $M_i$.  
Now, there exists a $m$ such that $\ap _k^m $ is trivial on 
$M_{n-1}$.  Suppose $(u_1, u_2,\cdots , u_k) \in V_{n-1}$.  Then $\ap _k^m 
(u_1, u_2,\cdots , u_k) = (u_1, u_2,\cdots , u_k)$.  This implies that 
for $1\leq i \leq k-1$, 
$u _i = u_i+\sum _{j>i}^k m_{i,j} u_j$ where $m_{i,j}\in {\mathbb N}$.  
For $i=k-1$, $u_{k-1} = u_{k-1} + m_{k-1, k}u_k$ and hence $u_k =0$.  
If $u_p =0$ for all $p>q>1$, then for $i=q-1$, $u_{q-1} = u_{q-1}+\sum 
_{j\geq q} m_{q-1, j}u_j = u_{q-1}+m_{q-1, q}u_q$ and hence $u_q =0$.  
Thus, $V_{n-1}$ is atmost one-dimensional.  Since $M_{n-1}/M_n = M_{n-1}$ 
is not finite, $M_{n-1}$ has dimension one and $V_{n-1} = \{ (u_1, 0, 
\cdots , 0) \mid u_1\in {\mathbb R} \}$.  

It can be seen that $T_k/M_{n-1} \simeq T_{k-1}$ and the action of 
${\mathbb Z}_{\ap _k}$ on $T_k/M_{n-1}$ is same as the action of ${\mathbb 
Z}_{\ap _{k-1}}$ on $T_{k-1}$.  Moreover, $M_{i+1}/M_{n-1} 
\subset M_i/M_{n-1}$ and ${M_i/M_{n-1}\over M_{i+1}/M_{n-1}} \simeq 
M_i/M_{i+1}$ for $0\leq i < n-1$ with $M_0/M_{n-1} = T_k /M_{n-1}$ and 
$M_{n-1}/M_{n-1} = (e)$.  By induction on $k$, we get that $n-1 = k-1$.  
\eo

Let $K=\prod _{k\in {\mathbb N}} T_k$.  Let $\ap \colon K \ra K$ be 
the automorphism defined by $\ap (f)(k) = \ap _k (f(k))$ for all $f 
\in K$ and all $k\in {\mathbb N}$.  Then $\ap$ is a continuous automorphism 
and the ${\mathbb Z}$-action defined by $\ap$ is distal on $K$.  

If there is a finite sequence $$(e) = K_{n} \subset K_{n-1} \subset \cdots 
\subset K_1 \subset K_0 =K$$ of $\ap$-invariant closed subgroups such that 
the action of ${\mathbb Z}_\ap$ on $K_i/K_{i+1}$ is finite for $i\geq 0$.  
This implies that each $T_k$ has a finite series 
$$(e) =K_{n,k}\subset K_{n-1, k}\subset \cdots 
\subset K_{1,k} \subset K_{0, k}= T_k$$ of 
$\ap _k$-invariant closed subgroups such that the action of 
${\mathbb Z}_{\ap _k}$ on $K_{i,k}/K_{i+1, k} $ is finite for $i\geq 0$.  

It follows from Lemma \ref{tl} that $k\leq n$.  Since $k\geq 1$ is 
arbitrary, this is a contradiction. Thus, the sequence $(K_i)$ of closed 
subgroups as in Proposition \ref{p1} for $K$ and $\ap$ is not finite.  

\end{section}

\begin{acknowledgement}
I would like to thank Prof. W. Jaworski for offering a Post-doctoral
fellowship of NSERC and many helpful suggestions and discussions. I would
also like to thank Department of Mathematics and Statistics, Carleton
University for their hospitality during my stay.
\end{acknowledgement}

\noindent {C. Robinson Edward Raja \newline
Stat-Math Unit \newline
Indian Statistical Instittue \newline
8th Mile Mysore Road \newline
Bangalore 560 059. India \newline
e-mail: creraja@isibang.ac.in}

\end{document}